\documentclass[MSNbibl,nameyear,dvips]{arxbj}

\aid{0}
\volume{19}
\issue{4}
\pubyear{2013}
\firstpage{1391}
\lastpage{1403}
\doi{10.3150/12-BEJSP07} 

\makeatletter

\newcommand{\cal}{\mathcal}

\newtheorem{algorithm}{Algorithm}

\makeatother

\begin{document}
\begin{frontmatter}

\title{Particle filters}
\runtitle{Particle filters}

\begin{aug}
\author{\fnms{Hans R.} \snm{K\"unsch}\corref{}\ead[label=e1]{kuensch@stat.math.ethz.ch}}
\runauthor{H.R. K\"unsch} 
\address{Seminar f\"ur Statistik, D-Math, ETH Zurich, CH-8092 Zurich,
Switzerland.\\
\printead{e1}}
\end{aug}


%
\begin{abstract}
This is a short review of Monte Carlo methods for
approximating filter distributions in state space models. The basic
algorithm and different strategies to reduce imbalance of the weights are
discussed. Finally, methods for more difficult problems like smoothing and
parameter estimation and applications outside the state space model context
are presented.
\end{abstract}

%

\end{frontmatter}

\section{Introduction}

Filtering is engineering terminology for extracting information about
a signal from partial and noisy observations. In geophysics, filtering
is usually called data assimilation. In the last 50 years,
filtering has been mainly studied in the framework of state space or hidden
Markov models, assuming a Markovian time evolution of the signal and
observations which are instantaneous functions of the signal subject to white
observation noise. Developments started in the 1960s with the Kalman--Bucy
filter (\citet{r25}, \citet{r26})
for linear Gaussian models and with the forward--backward
algorithm due to Baum and Welch for models with a finite state space
(see p. 74 of \citet{r9} for the history of this algorithm, including
references). The essential feature of these methods is that they
are recursive and thus suitable for online applications where the
observations arrive sequentially and quantities of interest have to be
recomputed with each new observation.

Probabilists started in the mid-sixties to develop a general theory of
nonlinear filtering in continuous time. In
statistics, state space models and filtering techniques took longer to
take roots. In the
seventies and eighties, the relation between linear state space and
ARMA models was studied and used. A breakthrough occurred
with the paper \citet{r23}
which developed recursive Monte Carlo methods called particle filters.
Interestingly,
\citet{r24} had proposed much earlier to use Monte Carlo methods,
but the idea of resampling was missing. However, this idea is essential
to ensure that the required sample size for a given accuracy does not explode
with the number of time steps. Particle filters quickly
became very popular. Among other things they have also been
used for continuous time filtering. Nowadays, they
are also applied outside the context of state space models as a complement
to other, static MCMC methods. In the 1990s,
geophysicists developed a different version of the particle filter,
called the Ensemble Kalman filter which is more stable in high dimensions.
After some delay, this idea
has now also become part of the research in statistics.

There are many presentations of the topic in books and in
survey articles (e.g., \citet{r27}, \citet{r16}, \citet{r13}, Capp{\'e}, Moulines and Ryd{\'e}n
(\citeyear{r9}), \citet{r8},
\citet{r12}, \citet{r17}). This paper gives a brief
introduction for
non-specialists, explaining the main algorithms, describing
their scope and also their limitations and surveying some of the
interesting current developments. Because of limitations of space,
many interesting topics and references that would deserve
to be mentioned had to be omitted.

\section{State space models}
\subsection{Definitions}
A state space model consists of an unobservable $(S, \mathcal{S})$-valued
Markov process $(X_t)$, the state of a system or the signal, combined
with partial and noisy $\mathbb{R}^d$-valued observations $(Y_i; i
\geq1)$ of the
state at
discrete times $t_i$. In order to simplify the notation, we assume
$t_i=i$. We denote the initial distribution of $X_0$ by
$\pi_0$ and the conditional distribution of $X_i$ given
$X_{i-1}=x_{i-1}$ by
$P(dx_i|x_{i-1})$. Observations at different
times are assumed to be conditionally independent given the states, and
the conditional distributions of $Y_i$ given $X_i$ are
assumed to have densities $g$ with respect to some reference measure
$\nu$
(usually the Lebesgue or the counting measure). Time homogeneity of
these conditional distributions is only assumed to simplify notation.

The state process
can be in continuous or discrete time. In the former case, the transition
kernel $P$ is usually not available analytically. For some of the
algorithms, this is not necessary, it is sufficient that we are able
to simulate from $P(\cdot|x)$ for any value $x$. Because some
applications have
a deterministic or partially deterministic state evolution, we do not
assume the existence of densities for $P$.

Throughout, notation like $X_{0:n}$ for $(X_0,\ldots, X_n)$ is used.
By a slight abuse of notation, $p$~stands for any (conditional)
density: The arguments of $p$ will indicate which random variables are involved.
The ratio of two probability measures is an abbreviation for the
Radon--Nikodym derivative.

\subsection{Examples}
State space models have a wide range of applications in finance
(stochastic volatility, interest rates), engineering (tracking, speech
recognition, computer vision), biology (genome sequence analysis, ion channels,
stochastic kinetic models), geophysics (meteorology, oceanography,
reservoir modeling), analysis of longitudinal data and others. It is not
possible here to describe these applications in detail or give references
to all relevant pulications.
Some of these applications are discussed in
\citet{r27} and in \citet{r16}. A~few references of more
recent applications are \citet{r6} and \citet{r32} for biology,
Part IX in \citet{r12} for financial mathematics,
and \citet{r19}, \citet{r1} and \citet{r10}
for geophysical applications.

\section{The basic particle filter}
\subsection{Filtering recursions}
By the assumptions on the state and observation process, we have the
following joint distributions for $n \geq m$
%
\begin{equation}
\label{eqjoint} (X_{0:n},Y_{1:m}) \sim\pi_0(dx_0)
\prod_{t=1}^n P(dx_t|x_{t-1})
\prod_{t=1}^m g(y_t|x_t)
\nu(dy_t).
\end{equation}
The information about the state contained in the observations
is expressed by the conditional distributions
$\pi_{s:t|n}$ of $X_{s:t}$ given $Y_{1:n}=y_{1:n}$. Of particular
interest are $\pi_{0:n}:=\pi_{0:n|n}$, called here the joint smoothing
distribution, and $\pi_{n}:=\pi_{n|n}$, called here the filter
distribution (the terminology is not unique). For $n \geq m$, $\pi_{0:n|m}$
follows immediately from (\ref{eqjoint}) and Bayes formula. Other cases
are then obtained
in principle by marginalization. We are however interested in methods
to compute or approximate expectations with respect to these
distributions in an explicit and efficient way. For this, recursive
formulae are
most useful. It is straightforward to verify that
%
\begin{eqnarray}
\label{eqprop1}
\pi_{0:n|n-1}(dx_{0:n}|y_{1:n-1}) &=& \pi_{0:n-1}(dx_{0:n-1}|y_{1:n-1})
P(dx_n|x_{n-1}) ,
\\
\label{equpdate1}
\pi_{0:n}(dx_{0:n}|y_{1:n}) &=& \pi_{0:n|n-1}(dx_{0:n}|y_{1:n-1})
\frac{g(y_n|x_n)}{p_n(y_n|y_{1:n-1})},
\end{eqnarray}
where
%
\begin{equation}
\label{eqlikeli} p_n(y_n|y_{1:n-1}) = \int
\pi_{n|n-1}(dx_{n}|y_{1:n-1}) g(y_n|x_n).
\end{equation}
By marginalization, we therefore also have the recursions
%
\begin{eqnarray}
\label{eqprop2}
\pi_{n|n-1}(dx_{n}|y_{1:n-1}) &=& \int
\pi_{n-1}(dx_{n-1}|y_{1:n-1}) P(dx_n|x_{n-1}),
\\
\label{equpdate2}
\pi_{n}(dx_{n}|y_{1:n}) &=& \pi_{n|n-1}(dx_{n}|y_{1:n-1})
\frac{g(y_n|x_n)}{p_n(y_n|y_{1:n-1})}.
\end{eqnarray}
In both cases, the recursions consist of
a propagation step (\ref{eqprop1}) or (\ref{eqprop2}), respectively,
and an update or correction step, (\ref{equpdate1}) or (\ref{equpdate2}),
respectively. Typically, one wants to compute these recursions for an
arbitrary, but fixed sequence $y_1,y_2,\ldots$ (not necessarily a
realization from the state space model).

\subsection{Analytical solutions}
There are two important special cases where one can perform the above
recursions exactly. In the first one, the state space $S$ is finite
and the integrals reduce to finite sums which can be computed with
$O(n |S|^2)$ operations.
The second special case are linear Gaussian state space models
where $X_n | X_{n-1} \sim{ \cal N}(F X_{n-1}, V)$ and
$Y_n | X_{n} \sim{\cal N}(H X_{n}, R)$. If $\pi_0$ is also Gaussian, then
all $\pi_n$ are Gaussian and (\ref{eqprop2})--(\ref{equpdate2})
lead to recursions for the
conditional means and covariances. For comparison with the Ensemble
Kalman filter below, we write down the update step for going
from $\pi_{n|n-1}= {\cal N}(m_{n|n-1},P_{n|n-1})$ to
$\pi_{n}= {\cal N}(m_{n},P_{n})$:
%
\begin{equation}
\label{eqkalman} m_n = m_{n|n-1} + K_n(Y_n
- H m_{n|n-1}), \quad P_n = P_{n|n-1} - K_n
H P_{n|n-1}
\end{equation}
where $K_n = P_{n|n-1}H'( H P_{n|n-1}H' + R)^{-1}$ is the so-called
Kalman gain.

In most other cases of practical interest, one has to approximate
the integrals involved in (\ref{eqprop2}) and (\ref{eqlikeli}).
Numerical approximations are difficult to use because the region
of main mass of $\pi_n$ changes with $n$ and is unknown in advance.
The particle filter tries to generate values in this region adaptively
as new observations arise.

\subsection{The particle filter}
The particle filter recursively computes importance sampling
approximations of $\pi_{n}$, that is
\[
\pi_{n}(dx_{n}|y_{1:n}) \approx\hat{
\pi}_{n}(dx_{n}|y_{1:n}) = \sum
_{i=1}^N W^i_n
\Delta_{X^i_n}(dx_n).
\]
Here the $W^i_n$ are random weights which sum to one, $X^i_n$ are random
variables called ``particles'' and $\Delta_x$ is the point mass at $x$.
At time $0$, we draw particles from $\pi_0$ and set $W^i_0=1/N$. At
time $n$ we start with $\hat{\pi}_{n-1}$ and draw independently new
particles $X^{i}_{n}$ from $P(\cdot|X^{i}_{n-1})$. By (\ref{eqprop2}),
the particles $X^{i}_{n}$ with weights $W^i_{n-1}$ provide an
importance sampling approximation of $\pi_{n|n-1}$. If we also
update the weights with $W^i_{n} \propto W^i_{n-1}g(y_{n}|X^i_{n})$,
we have closed the recursion by (\ref{equpdate2}).

This algorithm has the drawback that after a few iterations most particles
are located at positions very far away from the region of main mass of
$\pi_n$
and the weights are very unbalanced. This can be avoided by introducing
a resampling step before propagation such that particles with low weights
die and particles with high weights have much offspring that is
independently propagated afterwards. Thus the basic particle filter,
also called the bootstrap filter or SIR-filter (Sampling Importance
Resampling), works as follows.
%
\begin{algorithm}
1. Resample: Draw $(X^{*1}_{n-1},\ldots, X^{*N}_{n-1})$ from $\hat
{\pi}_{n-1}$.

2. Propagate: Draw $X^{i}_{n}$ from $P(\cdot|X^{*i}_{n-1})$,
independently for
different indices $i$.

3. Reweight: Set $W^i_{n} \propto g(y_{n}|X^i_{n})$.
\end{algorithm}
Note that for any function $\varphi: S \rightarrow\mathbb{R}$,
$N^{-1}\sum_i \varphi(X^{*i}_{n-1})$ always has a larger variance than
$\sum_i W^i_{n-1} \varphi(X^{i}_{n-1})$.
The advantage of resampling is seen only after one or several propagation
steps. Because of this, we resample at the beginning and not at the end
of a recursion.

As a byproduct, the particle filter gives also the following estimate
of (\ref{eqlikeli})
\[
\hat{p}_n(y_n|y_{1:n-1}) =\frac{1}{N} \sum
_{i=1}^N g\bigl(y_n|X^i_n
\bigr).
\]
One can show by induction that the
product $\prod_{t=1}^n \hat{p}_t(y_t|y_{1:t-1})$ is an exactly unbiased
estimator of $p(y_{1:n}) = \prod_{t=1}^n p_t(y_t|y_{1:t-1})$ for any
$n$ and any $y_{1:n}$, see Theorem 7.4.2 in \citet{r13}.
However, $\hat{p}_t$ is in general not unbiased
for $p_t$.

\subsection{Simple improvements}
Because the numbering of particles is irrelevant, we only need to know
the number of times $N^i_{n}$ that the $i$-th particle is selected
in the resampling step. One can therefore reduce the additional variability
introduced by resampling by a so-called balanced resampling scheme, meaning
that $\mathbb{E}(N^i_n) = N W^i_n$ and $|N^i_n - N W^i_n| < 1$.
The simplest such
scheme uses a uniform$(0,1)$ random variable $U$ and takes as
$N^i_n$ the number of points in the intersection of $(U + {\mathbb Z})/N$
with $(\sum_{k=1}^{i-1} W^k_n, \sum_{k=1}^{i} W^k_n]$. See
\citet{r11} for other
balanced resampling schemes. Since balanced resampling can always be used
at little extra cost, it is widely used.

A second improvement omits the resampling step whenever the
weights are sufficiently uniform. As criterion, one often uses
the so-called effective sample size which is defined as one
over $\sum_{i=1}^N (W^i_n)^2$,
see \citet{r30} for a justification of the name of this criterion.

In the propagation step, we can draw
$X^{i}_{n}$ not from $P(\cdot|X^{*i}_{n-1})$, but from any other
distribution $Q$ which dominates
$P(\cdot|X^{*i}_{n-1})$. We then have to adjust the weights in the
reweighting step. The correct weights are obtained
by setting $r \equiv1$ in step 4 of Algorithm~\ref{alg2} below.
By letting $Q$ depend not only
on $X^{*i}_{n-1}$, but also on the new observation $y_n$,
we can make the propagated particles $X^i_{n}$ more compatible with $y_{n}$
and thus the weights more balanced. In the so-called
auxiliary particle filter due to \citet{r31}, one uses the new observation
$y_{n}$ not only in the propagation step, but also in an additional reweighting
step before resampling. The goal of this additional reweighting is to bring
$\hat{\pi}_{n-1}$ closer to $\pi_{n-1|n}$.
Thus, the auxiliary particle filter works as follows.
%
\begin{algorithm}\label{alg2}
1. Reweight: Set
\[
\hat{\pi}_{n-1|n} = \sum_{i=1}^N W^{*i}_{n-1}
\Delta_{X^i_{n-1}}(dx_{n-1})
\]
where $W^{*i}_{n-1} \propto W^i_{n-1}
r(X^i_{n-1},y_{n})$.

2. Resample: Draw $(X^{*1}_{n-1},\ldots, X^{*N}_{n-1})$ from $\hat
{\pi}_{n-1|n}$.

3. Propagate: Draw $X^{i}_{n}$ from $Q(\cdot|X^{*i}_{n-1},y_{n})$,
independently for
different indices $i$.

4. Reweight: Set
\[
W^i_{n} \propto w^i_n:=
\frac{g(y_{n}|X^i_{n})} {
r(X^{*i}_{n-1},y_{n})} \frac
{P(dx_n|X^{*i}_{n-1})}{Q(dx_n|X^{*i}_{n-1},y_{n})}\bigl(X^i_n\bigr).
\]
\end{algorithm}
In order to understand the formula for $w^i_n$, note that
$(X^{*i}_{n-1},X^i_n)$ has distribution proportional to
$\hat{\pi}_{n-1}(dx_{n-1})r(x_{n-1},y_n)Q(dx_n|x_{n-1},y_{n})$ and the
distribution target is proportional to
$\hat{\pi}_{n-1}(dx_{n-1})g(y_n|\allowbreak x_n)P(dx_n|x_{n-1})$. Because the average
of the unnormalized weights $w^i_n$ estimates the ratio of the normalizing
constants, the estimate of (\ref{eqlikeli}) is now
\[
\hat{p}_n(y_n|y_{1:n-1}) = \frac{1}{N} \sum
_{i=1}^{N} w^i_n
\cdot\sum_{k=1}^{N} r\bigl(X^k_{n-1},y_n
\bigr) W^k_{n-1}.
\]
The product $\prod_{t=1}^n \hat{p}_t(y_t|y_{1:t})$ is again unbiased for
$p(y_{1:n})$.

Auxiliary particle filters cannot be used if the state evolution is
deterministic,
or if the density $dP(\cdot|x')/dQ(\cdot|x',y)$ is not available in
closed form.
In other cases, the choices of $r$ and $Q$ are up to the user. Ideally, we
take
$r(x,y) = \int g(y|x') P(dx'|x)$ and
$Q(dx'|x,y) = r(x,y)^{-1} g(y|x')P(dx'|x)$, because then the weights
$W^i_{n}$ in the fourth step are constant. In most cases, these
choices are not possible, but one can try to find suitable approximations.
With the ideal choices for $r$ and $Q$, the auxiliary particle filter therefore
leads to a reweighting with
$p(y_n|x_{n-1})$ instead of $p(y_n|x_{n})$: Although this usually
reduces the variance of the weights, the gain may not be substantial.
In principle, it is possible to go further back in time by computing
particle filter approximations of $\pi_{n-L:n|n}$ for some $L>0$.
An auxiliary particle
filter in this case uses $y_{n-L:n}$ to reweight the particles
at time $n-L-1$ and to generate new particles at times $n-L$ to $n$.

\section{Complications and solutions}
\subsection{Main difficulties}
The main difficulty with the particle filter is that often weights
become unbalanced, even when we use the auxiliary particle filter
in Algorithm~\ref{alg2} or
apply some of the other simple improvements discussed
above. In such cases, most resampled particles coincide
(``sample depletion''). If the state transitions are partially deterministic,
this becomes especially drastic because the propagation will not
create diversity.

Partially deterministic state transitions occur for instance if the
model contains unknown
parameters $\theta$ in the state transition $P$ or in the observation
density $g$ and one proceeds by considering the enlarged state vector
$(\theta, X_t)$. The propagation step for $\theta$ is then simply
$\theta^i_{n}=\theta^{*i}_{n-1}$. One can add some noise to create
diversity, possibly combined with some shrinking towards the mean
to keep the variance the same. Still, this does not always work well.

A second instance with partially deterministic state transitions
occurs if one uses the particle
filter algorithm to approximate not only $\pi_n$, but
the whole smoothing distribution~$\pi_{0:n}$.
In principle, this is straightforward: Each particle at time $n$ is
then a path of length $n+1$ that we write as $X^i_{0:n|n}$. The propagation
step concatenates
a resampled path $X^{*i}_{0:n-1|n-1}$ with a new value
$X^i_n \sim P(\cdot|X^{*i}_{n-1|n-1})$.

If the weights at one time point become very unbalanced,
the filter can be completely unreliable and it can lose track
even though the propagation step later creates again diversity.
Unbalanced weights have been observed to occur easily if the dimension of
the observations is large. A theoretical explanation of this phenomenon
has been provided by \citet{r5}.

In the following, we discuss some more advanced methods that have
been proposed to overcome these difficulties.

\subsection{Resample moves}
\citet{r22} have proposed the following method to avoid sample depletion
when the particle filter is used to
produce an approximation of $\pi_{0:n}$ with particles $X^i_{0:n|n}$
and equal weights. Let $K_n$ be a Markov kernel on $S^{n+1}$ which has
$\pi_{0:n}$ as invariant distribution, constructed for instance according
to the general Metropolis--Hastings recipe. Drawing new particles
$X^{*i}_{0:n|n} \sim K_n(.|X^i_{0:n|n})$,
independently for different $i$'s will then give a new approximation of
$\pi_{0:n}$ which is expected to be at least as good as the old one.
If $K_n$ modifies all components of $X^i_{0:n|n}$, this method
also removes ties, but since typically a single
kernel can only update one or a few components of $X^i_{0:n|n}$, the
computational complexity increases with $n$ if one wants to get rid
of all ties.

\subsection{Ensemble Kalman filter}
This method is due to \citet{r18}. It assumes linear observations
with Gaussian errors, that is, $g(\cdot|x)$ is a normal density with mean
$Hx$ and variance $R$. It uses particles with equal weights,
the propagation step is the same as in the particle filter whereas
the update step is a Monte Carlo implementation of the Kalman filter
update (\ref{eqkalman}) with estimated first and second moment of
$\pi_{n|n-1}$:
%
\begin{algorithm}
1. Propagate: Draw $X^{*i}_{n}$ from $P(\cdot|X^{i}_{n-1})$.

2. Update: Draw i.i.d. values $\varepsilon^i_{n} \sim{\cal N}(0,
R)$ and set
\[
X^i_{n} = X^{*i}_{n} +
\widehat{K}_{n}\bigl(y_{n} - HX^{*i}_{n}
+ \varepsilon^i_{n}\bigr)
\]
where $\widehat{K}_n$ is the Kalman gain computed with the sample
covariance $\widehat{P}_{n|n-1}$ of the $X^{*i}_{n}$'s.
\end{algorithm}
It is not difficult to show that the algorithm is consistent as $N
\rightarrow
\infty$ for a linear Gaussian state space model. However, for
non-Gaussian $\pi_{n|n-1}$, this update typically has a systematic error
because only the location, but neither the spread nor the shape of
the sample $(X^i_{n})$ change if $y_n$ changes.
Nevertheless, the
Ensemble Kalman filter is extremely wide-spread in geophysical applications
where the state evolution is usually complicated, making
the propagation step the computational bottleneck. This forces one to
use a sample size $N$ which is much smaller than the dimensions of the state
or the observation. Even in such cases, the Ensemble Kalman filter
turns out to be surprisingly robust -- provided we
regularize the estimate $\widehat{P}_{n|n-1}$ of the covariance of
$\pi_{n|n-1}$.

Several attempts have been made to find algorithms which combine the
robustness of the Ensemble Kalman filter with the nonparametric features
of the particle filter. They either approximate $\pi_{n|n-1}$ by a
mixture of
Gaussians or use the Ensemble Kalman filter as a proposal distribution $Q$
in a particle filter. See \citet{r21} for references
and a new proposal which avoids both the fitting of a Gaussian mixture
to the forecast sample $(X^{*i}_{n})$ and the estimation of the
density $dP/dQ$ (which is usually not known analytically in these
applications).

An extension of the Ensemble Kalman filter to more general observation
densities $g$ has been given in \citet{lei11}.

\subsection{Particle smoothing}
In an offline application where all $T$ observations are available
from the beginning, one
can use smoothing algorithms which combine a forward filtering pass
through the data from $n=0$ to $n=T$ with a backward recursion from
$n=T-1$ to $n=0$. We limit ourselves to approximations of the marginals
$\pi_{n|T}$, but the same methods apply also for joint distributions.

By Bayes formula and conditional independence, we obtain the following
relations
%
\begin{eqnarray}
\label{eqsmooth0}
\pi_{n|T}(dx_n|y_{1:T}) &=& \pi_{n|n-1}(dx_n|y_{1:n-1})
\frac{p(y_{n:T}|x_n)}{p(y_{n:T}|y_{1:n-1})}
\\
\label{eqsmooth1}
&=& \pi_n(dx_n|y_{1:n}) \frac{p(y_{n+1:T}|x_n)}{p(y_{n+1:T}|y_{1:n})}.
\end{eqnarray}
This is also called the two-filter formula because we have the recursions
%
\begin{eqnarray}
\label{eqprop-back}
p(y_{n+1:T}|x_n) &=& \int p(y_{n+1:T}|x_{n+1})
P(dx_{n+1}|x_n),
\\
\label{equpdate-back}
p(y_{n:T}|x_n) &=& g(y_n|x_n)
p(y_{n+1:T}|x_{n})
\end{eqnarray}
which are dual to (\ref{eqprop2}) and (\ref{equpdate2}). Combining
(\ref{eqsmooth0})--(\ref{eqsmooth1}) with (\ref{eqprop-back}) gives
%
\begin{equation}
\label{eqsmooth2} \frac{\pi_{n|T}(dx_n|y_{1:T})}{\pi_n(dx_n|y_{1:n})}
\propto\int\frac{\pi_{n+1|T}(dx_{n+1}|y_{1:T})}{\pi_{n+1|n}(dx_{n+1}|y_{1:n})}
P(dx_{n+1}|x_n).
\end{equation}

In order to be able to use Monte Carlo methods, we have to assume
that for any $x'$ the state transition kernel $P(\cdot|x')$
has density $p(\cdot|x')$ with respect to some measure $\mu$ on~$S$.
Then the filter\vadjust{\goodbreak} distributions also have densities which we denote
by the same symbol. The right-hand side of (\ref{eqsmooth2}) can then be
considered as an integral with respect to $\pi_{n+1|T}$. Thus we
obtain a marginal
particle smoother $\hat{\pi}_{n|T}$ which has the same particles
as the filter, but different weights $W^i_{n|T}$ which are computed
with the recursion
\[
W^i_{n|T} = \sum_{k=1}^N
W^k_{n+1|T} \frac{W^i_n p(X^k_{n+1}|X^i_n)} {
\sum_j W^j_n p(X^k_{n+1}|X^j_n)}.
\]
The disadvantage is the complexity of the algorithm
which is of the order $O(N^2)$.

The algorithm in \citet{r7} computes first backward particle
approximations of the
distributions $\bar{\pi}_n(dx_n|y_{n:T}) \propto
p(y_{n:T}|x_{n})h_n(x_n)\mu(dx_n)$
where $h_n$ is a known function such that $p(y_{n:T}|x_{n})h_n(x_n)$ is
integrable. Inserting a forward particle filter approximation for
$\pi_{n|n-1}$ and a backward particle filter approximation for
$p(y_{n:T}|x_{n})$ into (\ref{eqsmooth0}) gives then an approximation
of $\pi_{n|T}$ which is concentrated on the particles approximating
$\bar{\pi}_n$.

\citet{r20} have suggested to insert particle approximations
into
\[
\pi_{n|T}(x_n|y_{1:T}) \propto
\pi_{n|n-1}(x_n|y_{1:n-1}) g(x_n|y_n)
\int\frac{p(x_{n+1}|x_n)}{h_{n+1}(x_{n+1})} \bar{\pi}_{n+1}(dx_{n+1}|y_{n+1:T})
\]
which follows by combining (\ref{eqsmooth0}) with
(\ref{eqprop-back}) and (\ref{equpdate-back}). This has the
advantage that
the support of $\hat{\pi}_{n|T}$ is not constrained on the sampled
particles from the forward or the backward recursion. Moreover, one can
sample from the approximation with an algorithm of complexity $O(N)$
which may not be efficient, however.


\subsection{Particle MCMC}
This is a recent innovation by \citet{r3} which uses particle
filters as
a building block in an MCMC algorithm. Assume that $g$ and the density
of $P$ both depend on an unknown parameter $\theta$ with prior
density $p(\theta)$ and that we want
to sample from the posterior $p(x_{0:T},\theta|y_{1:T})$. A
Gibbs sampler which updates single components of $x_{0:T}$ given the
rest is usually too slow, and exact updates of the whole sequence
$x_{0:T}$ are usually not possible. What the particle filter provides
are \textit{random approximations} $\hat{\pi}_{0:T;\theta}$ of
$\pi_{0:T;\theta}=p(x_{0:T}|y_{1:T},\theta)$ for any fixed
$\theta$. \citet{r3} show that with these random approximations
one can still construct Markov chains which leave the
correct posterior invariant without letting the number of particles
go to infinity.

The first such algorithm is called particle marginal Metropolis--Hastings
sampler. It is an approximation of the sampler which jointly proposes
$(\theta',x'_{0:T})$ from the distribution
$q(\theta'|\theta) \,d\theta' \pi_{0:T}(dx'_{0:T}|y_{1:T},\theta')$
with the acceptance ratio
\[
\frac{p(y_{1:T}|\theta') p(\theta') q(\theta'|\theta)} {
p(y_{1:T}|\theta) p(\theta) q(\theta|\theta')}.
\]
The approximation occurs at two places: First $x'_{0:T}$ is generated
from $\hat{\pi}_{0:T;\theta'}$ instead of~$\pi_{0:T;\theta}$, and second
the unknown likelihoods\vadjust{\goodbreak}
$p(y_{1:T}|\theta')$ and $p(y_{1:T}|\theta)$ in the acceptance ratio
are replaced by unbiased estimates from the particle filter. The
surprising result is that the errors from these two approximations
cancel and the algorithm has the exact posterior
$p(x_{0:T},\theta|y_{1:T})$ as invariant distribution for any $N$.

Instead of jointly proposing a parameter and a path of the state process,
one can also use a Gibbs sampler, alternating between updates of the
parameter and the state process. Updating the parameter given the state
and the observations is usually feasible, but for the other update one
samples again from a particle filter approximation
$\hat{\pi}_{0:T;\theta}$ and not from $\pi_{0:T;\theta}$.
\citet{r3} show that this also gives a correct
algorithm for any $N > 1$ provided the particle filter approximation is
modified such that the current path is equal to one of the
particle paths $X^k_{0:T}$ in $\hat{\pi}_{0:T;\theta}$.

\section{Convergence results}
Laws of large numbers as well as central limit theorems have been shown
for particle filter approximations. \citet{r13} contains
general results, \citet{r28} gives an essentially self-contained
short derivation. First, one can show that
for every $n$, every $y_{1:n}$ and a suitable class of
functions $\varphi$, $\int\varphi(x) \hat{\pi}_n(dx_n|y_{1:n})$ converges
in probability or almost surely to $\int\varphi(x) \pi_n(dx_n|y_{1:n})$.
The proof works by induction on $n$, assuming that
$\hat{\pi}_{n-1}$ is close to $\pi_{n-1}$. This error propagates
in the next particle filter iteration, but one can control
by how much it grows in the worst case, and the additional Monte Carlo
error in the $n$-th step can be bounded by standard methods, at least
with multinomial (independent) resampling. For balanced sampling, there
seems to be still no general proof.

However, such a result is of limited use because the required sample
size $N$ may grow exponentially with the number of steps $n$.
For applications, it is more relevant to find conditions under which the
convergence is uniform in $n$. This is more difficult because -- in
contrast to the propagation step -- the update
step is in general not contractive and the above
induction argument does not succeed. One has instead to study the error
propagation over several time steps. This is equivalent to the question if
and how fast the filter forgets its initial distribution $\pi_0$ which has
been studied extensively, see e.g. \citet{r4}.

\section{More general situations}

\subsection{Filtering with continuous time observations}
Much of the probability literature on filtering considers
both state and observation processes in
continuous time. More precisely,
$(Y_t)$ is assumed to satisfy the following evolution equation
\[
dY_t = h(X_t)\,dt + dB_t
\]
where $(B_t)$ is a multivariate Brownian motion. We again denote by
$\pi_t$ the conditional distribution of $X_t$ given the
$\sigma$-field generated by the observations $(Y_s, s \leq t)$
(completed by all null sets).\vadjust{\goodbreak} Note that $(\pi_t)$ is a stochastic process
which takes values in the set of probability measures on
$(S,{\mathcal S})$. The evolution equation for $(\pi_t)$
corresponding to the recursions (\ref{eqprop2})--(\ref{equpdate2})
is a stochastic PDE, the Kushner--Stratonovich equation.
A particle filter approximation consists of interacting particles
$(X^i_t)$ and
associated weights $(W^i_t)$: Within an interval of length $\delta$
they evolve independently, whereas at multiples of $\delta$ there
is a resampling step like in the discrete case, see \citet{r33}
for more details.

\subsection{Sampling from moving targets}
Particle filtering algorithms have found many applications outside
the state space framework. In these cases, the more general term
sequential Monte Carlo is used. Assume we have a complicated target distribution
$\pi$ on $(S,{\mathcal S})$ which we cannot sample directly. In such a
situation,
a~promising strategy consists of sampling recursively from a sequence
$\pi_0, \pi_1,\ldots, \pi_T$ where $\pi_0$ is a simple distribution,
$\pi_T=\pi$ is the target one is interested in, and $\pi_n$ is close to
$\pi_{n-1}$. One example is the posterior distribution of a parameter
with a large number $T$ of observations where $\pi_n$ the
posterior for the first $n$ observations. In another example,
the $\pi_n$'s are tempered approximations of $\pi=\pi_T$:
\[
\pi_n(dx) \propto \biggl(\frac{\pi_T(dx)}{\pi_0(dx)} \biggr)^{\phi_n}
\pi_0(dx) \quad(0=\phi_0 < \phi_1 <
\cdots< \phi_T=1).
\]

Starting with a sample from $\pi_0$, one wants to recursively generate
samples $(X^i_n)$ from $\pi_n$ by resampling, propagation
and reweighting as in the particle filter. If $K_n$ denotes the transition
kernel in the $n$-th propagation step, then for reweighting we need the
density of $\pi_n$ with respect to $\int\pi_{n-1}(dx')K_n(\cdot
|x')$ which
is typically not available in closed form, unless we choose $K_n$
such that it leaves $\pi_{n-1}$ invariant. The idea in \citet{r14}
which allows more flexibility for the choice of $K_n$ is to
consider the distributions
$\pi_{n-1}(dx') K_n(dx|x')$ and $\pi_n(dx)L_n(dx'|x)$
on the product space $(S,{\mathcal S})^2$. Here
$L_n$ is an arbitrary kernel such that these two distributions are
absolutely continuous. If $(X^i_{n-1})$ is a (weighted) sample from
$\pi_{n-1}$
and we draw $X^i_n$ from $K_n(dx|X^i_{n-1})$ independently for different
$i$'s, then $(X^i_{n-1},X^i_{n})$ is a (weighted) sample from $\pi_{n-1}(dx')
K_n(dx|x')$. We can convert this into a weighted sample from
$\pi_n(dx)L_n(dx'|x)$ because the Radon--Nikodym density can be computed
without integration. By marginalization, we finally obtain the desired
weighted sample from $\pi_n$. In \citet{r14}, the optimal choice of
$L_n$ for given $K_n$ is determined.

\subsection{Rare event simulation}
Particle filters are also used in rare event simulation, see e.g.
\citet{r15}. Assume $(Z_t)$ is a Markov process with
fixed starting point $z_0$, $\tau$ and $\zeta$ are two stopping times
and we are interested in $\mathbb{P}(\tau< \zeta)$ which is small so simple
Monte Carlo is inefficient.
In a technique called ``importance splitting'' one introduces a
sequence of
stopping times $0 = \tau_0 < \tau_1 < \cdots< \tau_T=\tau$ and sets
$X_n = Z_{\tau_n}$. Moreover, we introduce ``observations''
$Y_n = 1_{[\tau_n < \zeta]}$.\vadjust{\goodbreak} Then for $y_1=y_2= \cdots=y_T =1$,
$\pi_n$
is the conditional distribution of $Z_{\tau_n}$ given $\tau_n < \zeta$,
and $p(y_{1:T})$ is the probability we would like to estimate. Hence, we
can estimate this probability unbiasedly with a particle filter since
it gives \textit{unbiased} estimates of $p(y_{1:n})$. Because the
observations are deterministic functions of the state, the resampling
step simply duplicates the particles with $\tau_n < \zeta$ until the
sample size is again $N$. \citet{r2} propose to control precision
instead of
computational effort. This means that in the $n$-th step we do
not propagate a fixed number of
particles and see how many of them satisfy $\tau_n < \zeta$, but rather
propagate particles until a fixed number of them satisfies
$\tau_n < \zeta$. One can still obtain an unbiased estimator, and in addition
this can increase the efficiency of the algorithm. Also in other
applications of the particle filter where all observations are
available from the beginning, it can be worthwile to aim for
a fixed precision instead of a fixed computational effort in each
iteration, using for instance accept-reject methods instead of importance
sampling (see \citet{r28}).


%

\printhistory

\end{document}